\documentclass[a4,12pt]{article}
\usepackage{amstext}
\usepackage{amssymb}
\title{Dynkin's Isomorphism with Sign Structure}
\author{Kshitij Khare}

\begin{document}
\maketitle

\begin{abstract}
The Dynkin isomorphism associates a Gaussian field to a Markov chain. These Gaussian fields 
can be used as priors for prediction and time series analysis. Dynkin's construction gives 
rise to Gaussian fields with all non-negative covariances. We extend Dynkin's construction 
(by introducing a sign structure on the Markov chain) to allow general covariance sign 
patterns. 
\end{abstract}

\section{Introduction}
Let $\mathbb{X} = \{X_t\}_{t \geq 0}$ be a reversible Markov process with a countable state  
space $\mathcal{X}$ and symmetric generator matrix $Q = ((q_{xy}))_{x,y \in\mathcal{X}}$, 
such that all states are transient. To ensure transience of $Q$, it is sufficient to assume 
that either 
\begin{equation} \label{transient1}
Q \mbox{ is irreducible and } \sum_{y \in \mathcal{X}} q_{xy} < 0 \mbox{ for atleast one } 
x \in \mathcal{X}. 
\end{equation}

\noindent
or 
\begin{equation} \label{transient2}
\sum_{y \in \mathcal{X}} q_{xy} < 0 \mbox{ for every } x \in \mathcal{X}. 
\end{equation}

\noindent
We provide a proof of this in Section $2$. Under both of these assumptions, $Q$ is 
non-conservative, i.e. from atleast one state there is a positive probability of going to 
an absorbing ``cemetery" $\Delta$ (not included in $\mathcal{X}$) and staying there 
forever. Dynkin \cite{dynkismph1} associated a Gaussian field $\{Z_x\}_{x \in \mathcal{X}}$ 
with variance-covariance matrix 
$$
\Sigma = -Q^{-1} 
$$

\noindent
with this Markov process and derived various intersesting properties of this 
correspondence. Since then, this correspondence has been used in several contexts (See 
Section $3$). All individual covariances of the Gaussian field $\{Z_x\}_{x \in 
\mathcal{X}}$ are non-negative in this construction. In this paper, we will extend Dynkin's 
construction to a larger class of variance-covariance matrices, which allow for positive as 
well as negative covariances. 

\indent
For this purpose, we introduce a ``sign-matrix" $\mathcal{S}$ such that 
$$
\mathcal{S} (x,x) = 1, \; \mathcal{S} (x,y) = \mathcal{S} (y,x), \; \mathcal{S} (x,y) \in 
\{-1, 1\} \; \forall x \neq y \in \mathcal{X}. 
$$

\noindent
If $\mathcal{S} (x,y) = 1$, the transition from $x$ to $y$ is called a positive 
transition. If $\mathcal{S} (x,y)= -1$, the transition from $x$ to $y$ is called a 
negative transition. Let $S_i$ denote the random time corresponding to the $i^{th}$ jump 
for the Markov process $\{X_t\}_{t \geq 0}, \; i = 1,2,3, ..........$. Define the 
``sign-process'' $\mathbb{H} = \{H_t\}_{t \geq 0}$ by 
$$
H_t = \prod_{i=1}^{\infty} \mathcal{S} (X_{S_{i-1}}, X_{S_i}) 1_{\{S_i \leq t\}} \; 
(\mbox{with } S_0 = 0 \mbox{ and } H_0 = 1). 
$$

\noindent
To describe in words, $H_t = 1$ if the number of negative transitions of $\mathbb{X}$ 
upto time $t$ is even and $H_t = -1$ if the number of negative transitions of 
$\mathbb{X}$ upto time $t$ is odd. Also, the transition to the ``cemetery" $\Delta$ from 
any state is a positive transition by default. 

\indent
Consider a Gaussian field $\{Z_x^\mathcal{S}\}_{x \in \mathcal{X}}$ with 
variance-covariance matrix 
\begin{equation} \label{dknsgvrcvn}
\Sigma^\mathcal{S} = (-Q \circ \mathcal{S})^{-1}. 
\end{equation}

\noindent
(We explain what we mean by $(-Q \circ \mathcal{S})^{-1}$ when $\mathcal{X}$ is countably 
infinite in Section $2$). As usual, let 
$$
l_t^x := \int_0^t 1_{\{X_s = x\}} ds, \; t \geq 0, x \in \mathcal{X} 
$$

\noindent
denote the occupation time of the Markov process $\{X_t\}_{t \geq 0}$ in the state $x$ till 
time $t$. We prove that for a realization of $\mathbb{X}$ independent of 
$\{Z_x^\mathcal{S}\}_{x \in \mathcal{X}}$ and for each bounded Borel measurable function 
$F: \mathbb{R}^{|\mathcal{X}|} \rightarrow \mathbb{R}$ and $x,y \in \mathcal{X}$, 
\begin{equation} \label{dknimphsgn}
{\bf{E}} \left[ Z_x^\mathcal{S} Z_y^\mathcal{S} F \left( \frac{(Z_u^\mathcal{S})^2}{2}, \; 
u \in \mathcal{X} \right) \right] = \int {\bf{E}} \left[ F \left( 
\frac{(Z_u^\mathcal{S})^2}{2} + l_\infty^u, \; u \in \mathcal{X} \right) H_\infty \right] 
d \mu_{xy}, 
\end{equation}

\noindent
where $\mu_{xy}$ is the conditional probability measure given that the process $\{X_t\}_{t 
\geq 0}$ enters the ``cemetery" $\Delta$ eventually with $y$ being the last state it stays 
in before being killed, scaled by a factor of $-Q^{-1} (x,y)$. We also prove identities for 
conditional prediction of the Gaussian field $\{Z_x^\mathcal{S}\}_{x \in \mathcal{X}}$ in 
terms of the Markov process $\mathbb{X}$. If $A \subset \mathcal{X}$ is finite, then 
\begin{equation} \label{dknpredexp}
{\bf{E}}[Z_b^\mathcal{S} \mid Z_a^\mathcal{S}, \; a \in A] = \sum_{a \in A} {\bf{E}}_b 
[1_{\{X_{R_A} = a\}} H_{R_A}] Z_a^\mathcal{S}, \; \forall b \in \mathcal{X} \setminus A, 
\end{equation}

\noindent
where $R_A$ is the first time (greater than or equal to $S_1$) when the Markov process 
$\mathbb{X}$ hits $A$, and 
\begin{equation} \label{dknpredcvn}
Cov(Z_b^\mathcal{S}, Z_{b'}^\mathcal{S} \mid Z_a^\mathcal{S}, \; a \in A) = {\bf{E}}_b 
\left[ \int_0^\infty 1_{\{X_s = b'\}} H_s 1_{\{s < R_A\}} \right], \; \forall b,b' \in 
\mathcal{X} \setminus A. 
\end{equation}

\noindent
Hence, the formulas for Gaussian field predictions in this case can be expressed elegantly 
in terms of quantities related to the Markov process $\mathbb{X}$. 

\indent
Note that, Dynkin's construction is a special case of (\ref{dknsgvrcvn}) with $\mathcal{S} 
(x,y) = 1, \; \forall x \neq y \in \mathcal{X}$. In this case $H_t = 1 \; \forall t \geq 
0$. Also, the version of (\ref{dknimphsgn}) when $\mathcal{S} (x,y) = 1 \; \forall x \neq y 
\in \mathcal{X}$ is known as the Dynkin's isomorphism theorem. It is remarkable that all 
changes that arise in the formulas as a result of introducing a ``sign-matrix" are 
reflected by just the ``sign-process" $\mathbb{H}$. Note that it is easy to keep track of 
$\{H_t\}_{t \geq 0}$ while simulating $\{X_t\}_{t \geq 0}$. 

\section{Preliminaries}
We clarify what we mean by inverse of an infinite matrix, atleast the ones that we are 
dealing with. Let $\tilde{Q}$ be an infinite matrix which can be written as $\tilde{R} (I - 
\tilde{P})$, where $\tilde{R}$ is a diagonal matrix with negative entries and $|\tilde{P}| 
:= ((|p_{ij}|))_{0 \leq i,j < \infty}$ is a sub-Markov matrix satisfying $\sum_{n=0}^\infty 
|\tilde{P}|^n < \infty$ (i.e. each entry of the matrix is finite). Then the matrix 
$(\sum_{n=0}^\infty \tilde{P}^n) \tilde{R}^{-1}$ satisfies 
$$
\tilde{Q} \left( \sum_{n=0}^{\infty} \tilde{P}^n \right) \tilde{R}^{-1} = \left( 
\sum_{n=0}^\infty \tilde{P}^n \right) \tilde{R}^{-1} \tilde{Q} = I. 
$$

\noindent
Hence, in such cases {\bf{we define}} 
$$
\tilde{Q}^{-1} = \left( \sum_{n=0}^\infty \tilde{P}^n \right) \tilde{R}^{-1}. 
$$

\indent
As in the introduction, let $\{X_t\}_{t \geq 0}$ be a reversible Markov process with a 
countable state space $\mathcal{X}$ with symmetric generator matrix $Q = ((q_{xy}))_{x,y 
\in \mathcal{X}}$ satisfying (\ref{transient1}) or (\ref{transient2}). Note that $q_{xy} = 
q_{yx} \geq 0$ for $x \neq y, \; q_{xx} < 0$ and $\sum_{y \in \mathcal{X}} q_{xy} \leq 0$. 
Let $\{Y_i\}_{i \geq 0} := \{X_{S_i}\}_{i \geq 0}$ be the embedded discrete-time Markov 
chain with one step transition probabilities 
$$
p_{xy} = - \frac{q_{xy}}{q_{xx}} \mbox{ for } x \neq y, \; p_{xx} = 0. 
$$

\noindent
Let $P := ((p_{xy}))_{x,y \in \mathcal{X}}$ and let $Q^{diag}$ denote the diagonal matrix 
with diagonal entries same as $Q$. Then 
\begin{equation} \label{cntrltndst}
Q = Q^{diag} (I - P). 
\end{equation}

\noindent
We describe a typical path of $\{X_t\}_{t \geq 0}$. The process starts at an initial state 
$Y_0$. The process stays at $Y_n$ during $[S_n, S_{n+1})$ for $n = 0,1,2, .......... $ and 
at a random time $\xi = S_{\eta + 1}$ jumps from the state $Y_{\eta}$ to the ``cemetery" 
$\Delta$ (not included in $\mathcal{X}$), and stays there forever. (Here $\xi$ takes 
non-negative real values and $\eta$ takes non-negative integer values). The value of $\xi$ 
(and hence $\eta$) can be infinite for certain sample paths, in which case the path does 
not terminate. Note that the difference $1 - \sum_{y \in \mathcal{X}} p_{xy}$ represents 
the probability $p(x, \Delta)$ of a jump from $x$ to the ``cemetery" $\Delta$ for the 
embedded Markov chain $\{Y_m\}_{m \geq 0}$. Also, conditional on $\{Y_m\}_{m \geq 0}$, the 
intermediate jump times $\{S_{i+1} - S_i\}_{i \geq 0}$ are independent and have 
$Exponential(-q_{Y_i Y_i})$ distribution for $i = 0,1,2, ..........$ 

\indent
Let us prove by the method of contradiction that under (\ref{transient1}) or 
(\ref{transient2}), $\sum_{n=0}^\infty P^n < \infty$ and hence $Q^{-1}$ exists. For this it 
is enough to show that all states are transient. Suppose $x \in \mathcal{X}$ is recurrent. 
The assumptions (\ref{transient1}) or (\ref{transient2}) imply that $\exists n$ such that 
$P^n (x, \Delta) > 0$. This implies by \cite[Theorem 3.4]{dutprtythy} that starting from 
the ``cemetery" $\Delta$ there is a positive probability of reaching the state $x$, which 
is a contradiction as the ``cemetery" $\Delta$ is absorbing. Hence, $\sum_{n=0}^\infty P^n 
< \infty$ and $Q^{-1} = (\sum_{n=0}^\infty P^n) (Q^{diag})^{-1}$ exists. 

\indent
The lemma below relates $Q^{-1}$ to the expected infinite occupation times for the process 
$\{X_t\}_{t \geq 0}$. 
\newtheorem{lem}{Lemma}[section]
\begin{lem}
\begin{eqnarray*}
-Q^{-1} (x,y) = {\bf{E}}_x [l_{\infty}^y], \; \forall x,y \in \mathcal{X}. 
\end{eqnarray*}
\end{lem}

\noindent
{\it Proof} Firstly, $Q = Q^{diag} (I - P)$ leads to 
$$
Q^{-1} = \left( \sum_{n=0}^\infty P^n \right) (Q^{diag})^{-1}. 
$$

\noindent
By decomposing the path of the Markov chain in terms of jump times, we get, 
\begin{eqnarray*}
{\bf{E}}_x \left[ \int_0^{\infty} 1_{\{X_s = y\}} ds \right] 
&=& {\bf{E}}_x \left[ \sum_{i=0}^{\infty} (S_{i+1} - S_i) 1_{\{X_{S_i} = y\}} \right]\\
&=& {\bf{E}}_x \left[ \sum_{i=0}^{\infty} (S_{i+1} - S_i) 1_{\{Y_i = y\}} \right]\\
&=& {\bf{E}}_x \left[ \sum_{i=0}^{\infty} {\bf{E}}_x \left[ (S_{i+1} - S_i) \mid 
\{Y_m\}_{m \geq 0} \right] 1_{\{Y_i = y\}} \right]\\
&=& {\bf{E}}_x \left[ \sum_{i=0}^{\infty} \frac{-1}{q_{Y_i Y_i}} 1_{\{Y_i = y\}} \right] 
\end{eqnarray*}

\noindent
The previous equality follows from the fact that conditional on $\{Y_m\}_{m \geq 0}$, 
the intermediate jump times $\{S_{i+1} - S_i\}_{i \geq 0}$ are independent and have 
$Exponential(-q_{Y_i Y_i})$ distribution for $i = 0,1,2, ..........$ 

\noindent
This gives
\begin{eqnarray*}
{\bf{E}}_x \left[ \int_0^{\infty} 1_{\{X_s = y\}} ds \right] 
&=& - \frac{1}{q_{yy}} {\bf{E}}_x \left[ \sum_{i=0}^{\infty} 1_{\{Y_i = y\}} \right]\\
&=& - \frac{1}{q_{yy}} \sum_{i=0}^{\infty} {\bf{P}}_x \{Y_i = y\}\\
&=& - \frac{1}{q_{yy}} \sum_{i=0}^{\infty} P^n (x,y) 
\end{eqnarray*}

\noindent
Note that ${\bf{P}}_x$ denotes the probability distribution starting at the state $x$, 
while $P$ is the transition matrix for the embedded Markov chain $\{Y_m\}_{m \geq 0}$. 
Hence, 
\begin{eqnarray*}
{\bf{E}}_x [l_{\infty}^y] 
&=& {\bf{E}}_x \left[ \int_0^{\infty} 1_{\{X_s = y\}} ds \right]\\
&=& -(Q^{diag} (I - P))^{-1} (x,y)\\
&=& -Q^{-1} (x,y) \mbox{       } \blacksquare
\end{eqnarray*}

\section{History}
In a series of papers Dynkin \cite{dynkismph1, dynkismph2, dynkismph3, dynkismph4} proposed 
and built on his construction as a connection between random fields and Markov processes. 
Marcus and Rosen have used the refined knowledge about Gaussian fields (eg. continuity of 
sample paths) to develop fine properties of symmetric Markov processes (eg. continuity of 
local times) using Dynkin's construction. Their book \cite{mrcsrsnmpl} gives a detailed and 
accessible account of their methods. Sheppard \cite{sheppardrk} uses Dynkin's isomorphism 
to give a proof of the Ray-Knight theorem on the Markovianity of one-dimensional 
diffusions. The properties of Markov processes can be utilized for analyzing the 
corresponding Gaussian fields. Ylvisaker \cite{yskrprtndg} uses Gaussian fields (amenable 
to Dynkin's isomorphism) as Bayesian priors for prediction and design problems, and makes 
use of the formulas relating the prediction properties of the Gaussian field to the 
corresponding Markov process. Bolthausen \cite{bolthausen} uses Dynkin's isomorphism as a 
tool in analyzing the limiting behaviour of the Gaussian free field. Eisenbaum 
\cite{eisenbaum2} and also Marcus and Rosen \cite{mrcsrsnmpl} have established variants of 
Dynkin's isomrphism. In the case of diffusions, Eisenbaum \cite{eisenbaum1} shows that 
Dynkin's isomorphism theorem and the Ray-Knight theorems can be derived from each other. 
In \cite{eikamarosh} the authors use an unconditional version of Dynkin's isomorphism to 
obtain a Ray-Knight theorem for a class of symmetric Markov processes. Diaconis and Evans 
\cite{dkndcnsens} proposed a different construction by looking at $-Q$ as the 
variance-covariance matrix instead of $-Q^{-1}$. Their construction yields Gaussian fields 
with negative individual covariances. 

\section{Generalization of Dynkin's Isomorphism}
We again consider a Markov process $\{X_t\}_{t \geq 0}$ with a countable state space 
$\mathcal{X}$ and with a generator matrix $Q$ as in Section $2$. We introduce a 
``sign-matrix" $\mathcal{S}$ such that 
$$
\mathcal{S} (x,x) = 1, \; \mathcal{S} (x,y) = \mathcal{S} (y,x), \; \mathcal{S} (x,y) \in 
\{-1,1\} \; \forall x \neq y \in \mathcal{X}. 
$$

\indent
As explained in the introduction, if $\mathcal{S} (x,y) = 1$, the transition from $x$ to 
$y$ is called a {\bf{positive transition}}. If $\mathcal{S} (x,y) = -1$, the transition 
from $x$ to $y$ is called a {\bf{negative transition}}. The ``sign-process" $\mathbb{H} = 
\{H_t\}_{t \geq 0}$ is defined by 
\begin{eqnarray*}
H_t = \prod_{i=1}^{\infty} \mathcal{S} (Y_{i-1}, Y_i) 1_{\{S_i \leq t\}}, \; (\mbox{with } 
H_0 = 1), 
\end{eqnarray*}

\noindent
and $H_t = 1$ if the number of negative transitions of $\mathbb{X}$ upto time $t$ is even 
and $H_t = -1$ if the number of negative transitions of $\mathbb{X}$ upto time $t$ is odd. 
Also, the transition to the ``cemetery" $\Delta$ from any state is a positive transition by 
default. 

\smallskip

\noindent
{\it Example} If $\mathcal{S} (x,y) = 1  \; \forall x \neq y \in \mathcal{X}$, then 
$H_t \equiv 1 \; \forall t \geq 0$ and $\Sigma = -Q^{-1}$. 

\smallskip

\noindent
{\it Example} If $\mathcal{S} (x,y) = -1 \; \forall x \neq y \in \mathcal{X}$, then 
$H_t = (-1)^{|\{i \geq 1: \; S_i \leq t\}|} \; \forall t \geq 0$ and $\Sigma = 
-(I + P)^{-1} (Q^{diag})^{-1} \; (\because \mbox{ From } (\ref{cntrltndst}))$. 
For this particular case, the hyper-process $\{H_t\}_{t \geq 0}$ is $1$ between 
$[S_{2i}, S_{2i+1})$ and $-1$ between $[S_{2i+1}, S_{2i+2})$ for $i \geq 0$. 

\smallskip

\noindent
Define 
\begin{eqnarray*}
\tilde{l}_t^x := \int_0^t 1_{\{X_s = x, H_s = 1\}} ds - \int_0^t 1_{\{X_s = x, H_s = -1\}} 
ds = \int_0^t 1_{\{X_s = 1\}} H_s ds, \; \; t \geq 0, \; x \in \mathcal{X}. 
\end{eqnarray*}

\noindent
We interpret $\tilde{l}_t^x$ as the {\bf{net occupation time}} in the state $x$ till 
time $t$ (with $\int_0^t 1_{\{X_s = x, H_s = 1\}} ds$ and 
$\int_0^t 1_{\{X_s = x, H_s = -1\}} ds$ interpreted as the 
negative and positive occupation times respectively).  

\noindent
Define the matrix 
$$
\Sigma^\mathcal{S} := (-Q \circ \mathcal{S})^{-1}. 
$$ 

\noindent
Here $\circ$ denotes Hadamard product i.e. elementwise product of the two matrices. Note 
that $|P \circ \mathcal{S}| = P$ and $\sum_{n=0}^\infty P^n < \infty$. Hence $(-Q \circ 
\mathcal{S})^{-1}$ exists. Also, since $-Q \circ \mathcal{S}$ is a diagonally dominant 
matrix with positive diagonal entries, hence $\Sigma^\mathcal{S}$ is positive definite. 
Note that, an infinite matrix is defined to be positive definite if all its finite 
principal submatrices are positive definite. 

\begin{lem} \label{vrcvnintpn}
\begin{eqnarray*}
\Sigma^\mathcal{S} (x,y) = {\bf{E}}_x \left[ \int_0^{\infty} 1_{\{X_s = y\}} H_s ds \right] 
= {\bf{E}}_x [\tilde{l}_{\infty}^y], \; \forall x,y \in \mathcal{X}. 
\end{eqnarray*}

\noindent
i.e. $\Sigma^\mathcal{S} (x,y)$ is the expected net occupation time at $y$ starting at $x$. 
\end{lem}

\noindent
{\it Proof} Firstly, 
\begin{eqnarray*}
Q \circ \mathcal{S} 
&=& (Q^{diag} (I - P)) \circ \mathcal{S}\\
&=& Q^{diag} (I - P \circ \mathcal{S}) 
\end{eqnarray*}

\noindent
This gives 
$$
\Sigma^\mathcal{S} = (-Q \circ \mathcal{S})^{-1} = (P \circ \mathcal{S} - I)^{-1} 
(Q^{diag})^{-1}. 
$$

\noindent
By decomposing the path of the Markov chain in terms of jump times, we get, 
\begin{eqnarray*}
{\bf{E}}_x \left[ \int_0^{\infty} 1_{\{X_s = y\}} H_s ds \right] 
&=& {\bf{E}}_x \left[ \sum_{i=0}^{\infty} (S_{i+1} - S_i) 1_{\{X_{S_i} = y\}} 
H_{S_i} \right]\\
&=& {\bf{E}}_x \left[ \sum_{i=0}^{\infty} (S_{i+1} - S_i) 1_{\{Y_i = y\}} 
H_{S_i} \right]\\
&=& {\bf{E}}_x \left[ \sum_{i=0}^{\infty} {\bf{E}}_x \left[ (S_{i+1} - S_i) 
\mid \{Y_m\}_{m \geq 0} \right] 1_{\{Y_i = y\}} H_{S_i} \right] 
\end{eqnarray*}

\noindent
The previous equality follows from the fact that 
$$
H_{S_i} = \prod_{j=1}^{i} \mathcal{S} (Y_{j-1}, Y_j) 
$$

\noindent
is a function of $\{Y_m\}_{m \geq 0}$. This gives, 
\begin{eqnarray*}
{\bf{E}}_x \left[ \int_0^{\infty} 1_{\{X_s = y\}} H_s ds \right] 
&=& {\bf{E}}_x \left[ \sum_{i=0}^{\infty} \frac{-1}{q_{Y_i Y_i}} 1_{\{Y_i = y\}} 
H_{S_i} \right]\\
&=& - \frac{1}{q_{yy}} {\bf{E}}_x \left[ \sum_{i=0}^{\infty} 1_{\{Y_i = y\}} 
H_{S_i} \right]\\
&=& - \frac{1}{q_{yy}} \sum_{i=0}^{\infty} {\bf{E}}_x \left[ 1_{\{Y_i = y\}} 
H_{S_i} \right] 
\end{eqnarray*}

\noindent
The exchange of sum and expectation is justified by the fact that $|H_s| = 1$ and 
$\sum_{i=0}^{\infty} {\bf{E}}_x \left[ 1_{\{Y_i = y\}} \right] = (I - P)^{-1} 
(x,y) < \infty$. 

\noindent
Let us calculate ${\bf{E}}_x \left[ 1_{\{Y_i = y\}} H_{S_i} \right]$. 
\begin{eqnarray*}
{\bf{E}}_x \left[ 1_{\{Y_i = y\}} H_{S_i} \right] 
&=& {\bf{E}}_x \left[ 1_{\{Y_i = y\}} \prod_{j=1}^i \mathcal{S} (Y_{j-1}, Y_j) 
\right]\\
&=& \sum_{y_1, y_2, ... , y_{i-1} \in \mathcal{X}} \prod_{j=1}^i p_{y_{j-1} y_j} 
\prod_{j=1}^i \mathcal{S} (y_{j-1}, y_j) \mbox{ where } y_0 = x, y_i = y.\\
&=& \sum_{y_1, y_2, ... , y_{i-1} \in \mathcal{X}} \prod_{j=1}^i p_{y_{j-1} y_j} 
\mathcal{S} (y_{j-1}, y_j) \mbox{ where } y_0 = x, y_i = y.\\
&=& (P \circ \mathcal{S})^i (x,y) 
\end{eqnarray*}

\noindent
Hence, 
\begin{eqnarray*}
{\bf{E}}_x [{\tilde{l}}_{\infty}^y] 
&=& {\bf{E}}_x [\int_0^{\infty} 1_{\{X_s = y\}} H_s ds]\\
&=& - \frac{1}{q_{yy}} \sum_{i=0}^{\infty} (P \circ \mathcal{S})^i (x,y)\\
&=& - \frac{1}{q_{yy}} (I - P \circ \mathcal{S})^{-1} (x,y)\\
&=& -(Q^{diag} (I - P \circ \mathcal{S}))^{-1} (x,y)\\
&=& \Sigma^\mathcal{S} (x,y) 
\end{eqnarray*}

\noindent
The proof is complete. \hspace{20mm} $\blacksquare$ 

\smallskip

\noindent
We next prove the isomorphism theorem (\ref{dknimphsgn}) for a zero mean Gaussian process 
$\{Z_x^\mathcal{S}\}_{x \in \mathcal{X}}$ with variance-covriance matrix $\Sigma = (-Q 
\circ \mathcal{S})^{-1}$ and an independent realization $\{X_t\}_{t \geq 0}$ of the Markov 
process with generator $Q$. 

\subsection{The finite case}
We consider the case when $\mathcal{X}$ is finite. We proceed similarily as Dynkin 
\cite{dynkismph1} and first consider functions of the form 
\begin{eqnarray*}
F_{\underline{d}} (\underline{w}) = e^{- \sum_{u \in \mathcal{X}} d_u w_u}, 
\end{eqnarray*}

\noindent
where $\underline{d} = \{d_u\}_{u \in \mathcal{X}}$ is arbitrary with $d_u \geq 0, \; 
\forall u \in \mathcal{X}$. Let $D$ denote the diagonal matrix with diagonal entries 
$\{d_u\}_{u \in \mathcal{X}}$. Then, 
\begin{eqnarray*}
{\bf{E}} \left[ Z_x Z_y e^{- \sum_{u \in \mathcal{X}} d_u \frac{Z_u^2}{2}} \right] 
&=& \int_{{\mathbb{R}}^{|\mathcal{X}|}} \frac{z_x z_y}{(\sqrt{2 \pi})^{|\mathcal{X}|} 
det(-Q \circ \mathcal{S})} e^{-\frac{\underline{z}^T D \underline{z} - 
\underline{z}^T (Q \circ \mathcal{S}) \underline{z}}{2}} d \underline{z}\\
&=& \frac{det(D - Q \circ \mathcal{S})}{det(-Q \circ \mathcal{S})} 
\int_{{\mathbb{R}}^{|\mathcal{X}|}} \frac{z_x z_y}{(\sqrt{2 \pi})^{|\mathcal{X}|} 
det(D - Q \circ \mathcal{S})} e^{-\frac{\underline{z}^T (D - Q \circ \mathcal{S}) 
\underline{z}}{2}} d \underline{z}\\
&=& \frac{det((D - Q) \circ \mathcal{S})}{det(-Q \circ \mathcal{S})} 
\int_{{\mathbb{R}}^{|\mathcal{X}|}} \frac{z_x z_y}{(\sqrt{2 \pi})^{|\mathcal{X}|} 
det((D - Q) \circ \mathcal{S})} e^{-\frac{\underline{z}^T 
((D - Q) \circ \mathcal{S}) \underline{z}}{2}} d \underline{z}\\
&=& \frac{det((D - Q) \circ \mathcal{S})}{det(-Q \circ \mathcal{S})} 
((D - Q) \circ \mathcal{S})^{-1} (x,y) 
\end{eqnarray*}

\noindent
Note that ${\bf{E}} \left[ e^{-\sum_{u \in \mathcal{X}} d_u \frac{Z_u^2}{2}} \right] = 
\frac{det((D - Q) \circ \mathcal{S})}{det(-Q \circ \mathcal{S})}$ by a similiar 
calculation as above. Hence, 
\begin{equation} \label{gaussinteg}
{\bf{E}} \left[ Z_x Z_y e^{- \sum_{u \in \mathcal{X}} d_u \frac{Z_u^2}{2}} \right] = 
{\bf{E}} \left[ e^{- \sum_{u \in \mathcal{X}} d_u \frac{Z_u^2}{2}} \right] 
((D - Q) \circ \mathcal{S})^{-1} (x,y) 
\end{equation}

\noindent
Note that $-(D - Q)$ is the generator of a Markov process $\{\bar{X}_t\}_{t \geq 0}$ 
with the same structure as $\{X_t\}_{t \geq 0}$ except that at every state $x \in 
\mathcal{X}$, there is an additional killing rate of $d_x$. Let 
$\{\bar{Y}_m\}_{m \geq 0}$ be the embedded discrete-time Markov chain and 
$\{\bar{H}_t\}_{t \geq 0}$ the hyper-process corresponding to 
$(\{\bar{X}_t\}_{t \geq 0}, \mathcal{S})$. 

\noindent
Let us establish the change of measure formula from $\{\bar{Y}_m\}_{m \geq 0}$ to 
$\{Y_m\}_{m \geq 0}$. 
\begin{eqnarray*}
P_x \{ \bar{Y}_1 = y_1, \bar{Y}_2 = y_2, ... , \bar{Y}_n = y_n \} 
&=& \prod_{i=1}^n \frac{q_{y_{i-1} y_i}}{-q_{y_{i-1} y_{i-1}} + d_{y_{i-1}}} \; 
(\mbox{with } y_0 = x)\\
&=& \left( \prod_{i=1}^n \frac{-q_{y_{i-1} y_{i-1}}}{-q_{y_{i-1} y_{i-1}} + 
d_{y_{i-1}}} \right) P_x \{ Y_1 = y_1, Y_2 = y_2, ... , Y_n = y_n \} 
\end{eqnarray*}

\noindent
Hence, 
\begin{equation} \label{measurechg}
{\bf{E}}_x [F(\bar{Y}_1, \bar{Y}_2, ... , \bar{Y}_n)] = 
{\bf{E}}_x \left[ \left( \prod_{i=1}^n 
\frac{-q_{y_{i-1} y_{i-1}}}{-q_{y_{i-1} y_{i-1}} + d_{y_{i-1}}} \right) 
F(Y_1, Y_2, ... , Y_n) \right] 
\end{equation}

\noindent
for each bounded Borel measurable function $F$. 

\medskip

\noindent
As in the introduction, let us define the measure $\mu_{xy}$ by 
\begin{equation} \label{dftncndmsr}
\mu_{xy} \{C\} = -Q^{-1} (x,y) {\bf{P}}_x \{C \mid \eta < \infty, Y_{\eta} = y\}. 
\end{equation}

\noindent
It is the appropriately scaled conditional probability measure given that the process 
$\mathbb{X}$ enters the ``cemetery" $\Delta$ eventually with $y$ being the last state it 
stays in before being killed. Note that, 
\begin{eqnarray*}
& & \mu_{xy} \{ Y_1 = y_1, Y_2 = y_2, Y_3 = y_3, ... , Y_n = y_n, \eta = n \}\\
&=& -Q^{-1} (x,y) \frac{P_x \{ Y_1 = y_1, Y_2 = y_2, Y_3 = y_3, ... , Y_n = y_n, 
Y_{n+1} = \Delta \} 1_{\{y_n = y\}}}{P_x \{ \eta < \infty, Y_{\eta} = y \}}\\
&=& -Q^{-1} (x,y) \frac{P_x \{ Y_1 = y_1, Y_2 = y_2, ... , Y_n = y_n \} 
1_{\{y_n = y\}} p(y, \Delta)}{\sum_{n=0}^{\infty} 
P_x \{ Y_n = y, Y_{n+1} = \Delta \}}\\
&=& \frac{-Q^{-1} (x,y) p(y, \Delta)}{\sum_{n=0}^{\infty} P_x \{ Y_n = y \} 
p(y, \Delta)} P_x \{ Y_1 = y_1, Y_2 = y_2, ... , Y_n = y_n \} 1_{\{y_n = y\}}\\
&=& \frac{1}{-q_{yy}} P_x \{ Y_1 = y_1, Y_2 = y_2, ... , Y_n = y_n \} 
1_{\{y_n = y\}} 
\end{eqnarray*}

\noindent
The previous equality follows from the fact that $-Q^{-1} (x,y) = 
\frac{1}{-q_{yy}} \sum_{n=0}^{\infty} P_x \{ Y_n = y \}$. 

\smallskip

\noindent
Hence, 
\begin{equation} \label{rltnmrkvmu}
\int F(Y_1, Y_2, ... , Y_{\eta}) 1_{\{\eta = n\}} d \mu_{xy} = 
\frac{1}{-q_{yy}} {\bf{E}}_x \left[ F(Y_1, Y_2, ... , Y_n) 
1_{\{Y_n = y\}} \right] 
\end{equation}

\noindent
for each bounded Borel measurable function $F$. 

\medskip

\noindent
We combine these results to evaluate $((D - Q) \circ \mathcal{S})^{-1} (x,y)$ in terms of 
$\{l_{\infty}^u\}_{u \in \mathcal{X}}$ and $H_{\infty}$. Let $\mu_{xy}$ be the measure 
defined in (\ref{dftncndmsr}). 
\begin{lem} \label{lpctmidnty}
\begin{eqnarray*}
((D - Q) \circ \mathcal{S})^{-1} (x,y) = \int e^{- \sum_{u \in \mathcal{X}} d_u 
l_{\infty}^u} H_{\infty} d \mu_{xy}. 
\end{eqnarray*}
\end{lem}

\noindent
{\it Proof} Firstly, we observe that under $\mu_{xy}, \; \eta < \infty$ and hence 
$H_{\infty} = H_{S_{\eta}}$, which is a measurable function of 
$Y_1, Y_2, ... , Y_{\eta}$. Also, 
\begin{eqnarray*}
\sum_{u \in \mathcal{X}} d_u l_{\infty}^u = \sum_{i=0}^{\eta} d_{Y_i} 
(S_{i+1} - S_i). 
\end{eqnarray*}

\noindent
Hence, 
\begin{eqnarray*}
& & \int e^{- \sum_{u \in \mathcal{X}} d_u l_{\infty}^u} H_{\infty} d \mu_{xy}\\
&=& \sum_{n=0}^{\infty} \int e^{- \sum_{i=0}^n d_{Y_i} (S_{i+1} - S_i)} H_{S_n} 
1_{\{\eta = n\}} d \mu_{xy}\\
&=& \sum_{n=0}^{\infty} \frac{1}{-q_{yy}} {\bf{E}}_x \left[ e^{- \sum_{i=0}^n 
d_{Y_i} (S_{i+1} - S_i)} H_{S_n} 1_{\{Y_n = y\}} \right] \; (\because 
\mbox{From (\ref{rltnmrkvmu})})\\
&=& \sum_{n=0}^{\infty} \frac{1}{-q_{yy}} {\bf{E}}_x \left[ {\bf{E}}_x \left[ 
e^{- \sum_{i=0}^n d_{Y_i} (S_{i+1} - S_i)} \mid \{Y_m\}_{m \geq 0} \right] 
H_{S_n} 1_{\{Y_n = y\}} \right]\\
&=& \sum_{n=0}^{\infty} \frac{1}{-q_{yy}} {\bf{E}}_x \left[ \left( 
\prod_{i=0}^n \frac{-q_{Y_i Y_i}}{-q_{Y_i Y_i} + d_{Y_i}} \right) H_{S_n} 
1_{\{Y_n = y\}} \right] 
\end{eqnarray*}

\noindent
The previous equality follows from the fact that conditioned on $\{Y_m\}_{m \geq 0}$, 
the intermediate jump times $\{S_{i+1} - S_i\}_{i \geq 0}$ are independent 
and have $Exponential(-q_{Y_i Y_i})$ distribution for 
$i = 0,1,2, .......... $. 

\noindent
Hence, with $\{\bar{S}_i\}_{i \geq 0}$ denoting the random transition times for 
$\{\bar{X}_t\}_{t \geq 0}$, we get, 
\begin{eqnarray*}
\int e^{- \sum_{u \in \mathcal{X}} d_u l_{\infty}^u} H_{\infty} d \mu_{xy} 
&=& \sum_{n=0}^{\infty} \frac{1}{-q_{yy} + d_y} {\bf{E}}_x \left[ \left( \prod_{i=1}^n 
\frac{-q_{Y_{i-1} Y_{i-1}}}{-q_{Y_{i-1} Y_{i-1}} + d_{Y_{i-1}}} \right) H_{S_n} 
1_{\{Y_n = y\}} \right]\\
&=& \sum_{n=0}^{\infty} \frac{1}{-q_{yy} + d_y} {\bf{E}}_x \left[ \bar{H}_{\bar{S}_n} 
1_{\{\bar{Y}_n = y\}} \right] \; ( \because \mbox{From (\ref{measurechg})})\\
&=& {\bf{E}}_x \left[ \sum_{n=0}^{\infty} \frac{1_{\{\bar{Y}_n = y\}} 
\bar{H}_{\bar{S}_n}}{-q_{yy} + d_y} \right]\\
&=& {\bf{E}}_x \left[ \int_0^{\infty} 1_{\{\bar{X}_s = y\}} \bar{H}_s ds \right]\\
&=& ((D - Q) \circ \mathcal{S})^{-1} (x,y) 
\end{eqnarray*}

\noindent
The previous equality follows from the fact that $\{\bar{X}_t\}_{t \geq 0}$ is a Markov 
process with generator $-(D - Q)$ which satisfies (\ref{transient1}) or (\ref{transient2}). 

\noindent
Hence proved. \hspace{20mm} $\blacksquare$ 

\medskip

\noindent
It follows from this claim and (\ref{gaussinteg}) that 
\begin{eqnarray*}
{\bf{E}} \left[ Z_x^\mathcal{S} Z_y^\mathcal{S} e^{- \sum_{u \in \mathcal{X}} d_u 
\frac{(Z_u^\mathcal{S})^2}{2}} \right] = {\bf{E}} \left[ e^{- \sum_{u \in \mathcal{X}} d_u 
\frac{(Z_u^\mathcal{S})^2}{2}} \right] \int e^{- \sum_{u \in \mathcal{X}} d_u l_{\infty}^u} 
H_{\infty} d \mu_{xy}. 
\end{eqnarray*}

\noindent
Hence, 
\begin{eqnarray*}
{\bf{E}} \left[ Z_x^\mathcal{S} Z_y^\mathcal{S} e^{- \sum_{u \in \mathcal{X}} d_u 
\frac{(Z_u^\mathcal{S})^2}{2}} \right] = \int {\bf{E}} \left[ e^{- \sum_{u \in \mathcal{X}} 
d_u \left( \frac{(Z_u^\mathcal{S})^2}{2} + l_{\infty}^u \right)} H_{\infty} \right] d 
\mu_{xy}. 
\end{eqnarray*}

\noindent
The set of functions $F_{\underline{d}} (\underline{w}) = e^{- \sum_{u \in \mathcal{X}} 
d_u w_u}$, where $\underline{d} = \{d_u\}_{u \in \mathcal{X}}$ is arbitrary with 
$d_u \geq 0, \; \forall u \in \mathcal{X}$, generate the Borel $\sigma$-algebra in 
${\mathbb{R}}^{|\mathcal{X}|}$ and they form a closed class under multiplication. Also, 
the set of functions $F$ for which 
\begin{eqnarray*}
{\bf{E}} \left[ Z_x^\mathcal{S} Z_y^\mathcal{S} F \left( \frac{(Z_u^\mathcal{S})^2}{2}, \; 
u \in \mathcal{X} \right) \right] = \int {\bf{E}} \left[ F \left( 
\frac{(Z_u^\mathcal{S})^2}{2} + l_{\infty}^u, \; u \in \mathcal{X} \right) H_{\infty} 
\right] d \mu_{xy}, 
\end{eqnarray*}

\noindent
is a linear space closed under bounded convergence and under monotone convergence. 
Hence, for each bounded Borel measurable function $F: {\mathbb{R}}^{|\mathcal{X}|} 
\rightarrow \mathbb{R}$, 
\begin{equation} \label{dknismphss}
{\bf{E}} \left[ Z_x^\mathcal{S} Z_y^\mathcal{S} F \left( \frac{(Z_u^\mathcal{S})^2}{2}, \; 
u \in \mathcal{X} \right) \right] = \int {\bf{E}} \left[ F \left( 
\frac{(Z_u^\mathcal{S})^2}{2} + l_{\infty}^u, \; u \in \mathcal{X} \right) H_{\infty} 
\right] d \mu_{xy}, 
\end{equation}

\noindent
where $\{Z_x^\mathcal{S}\}_{x \in \mathcal{X}}$ is a zero mean Gaussian field with 
variance-covariance matrix given by $(-Q \circ \mathcal{S})^{-1}$, and 
$\{l_{\infty}^u\}_{u \in \mathcal{X}}$ are the occupation times of a realization 
$\{X_t\}_{t \geq 0}$ independent of $\{Z_x^\mathcal{S}\}_{x \in \mathcal{X}}$ of a 
Markov process with generator $Q$. Also, in this case the map $(Q, \mathcal{S}) \rightarrow 
(-Q \circ \mathcal{S})^{-1}$ is one-to-one, because 
\begin{eqnarray*}
(-Q_1 \circ {\mathcal{S}}_1)^{-1} = (-Q_2 \circ {\mathcal{S}}_2)^{-1} 
&\Leftrightarrow& Q_1 \circ {\mathcal{S}}_1 = Q_2 \circ {\mathcal{S}}_2\\
&\Leftrightarrow& Q_1 = Q_2, \; {\mathcal{S}}_1 = {\mathcal{S}}_2 
\end{eqnarray*}

\noindent
The previous statement is justified by the fact that $Q_i, \; i = 1,2$ have 
negative off-diagonal entries and ${\mathcal{S}}_i, \; i = 1,2$ have entries 
that equal $1$ or $-1$. 

\medskip

\noindent
Let us now turn our attention to the problem of predicting the above Gaussian field given 
observations in a proper subset $A \subset \mathcal{X}$, and proving the identities 
(\ref{dknpredexp}) and (\ref{dknpredcvn}). We do not require the assumption of independence 
of $\{X_t\}_{t \geq 0}$ and $\{Z_x^\mathcal{S}\}_{x \in \mathcal{X}}$ for these 
calculations. Let $B := \mathcal{X} \setminus A$. Note that, 
$$
{\bf{E}}[Z_B^\mathcal{S} \mid Z_A^\mathcal{S}] = {\Sigma}_{BA} {{\Sigma}_{AA}}^{-1} 
Z_A^\mathcal{S}. 
$$

\noindent
Since $\Sigma^\mathcal{S} = (-Q \circ \mathcal{S})^{-1}$, it follows that 
$$
{\Sigma}_{BA}^\mathcal{S} {{\Sigma}_{AA}^\mathcal{S}}^{-1} = -\{(-Q \circ
\mathcal{S})_{BB}\}^{-1} (-Q \circ \mathcal{S})_{BA}. 
$$

\noindent
By slightly detailed, but straightforward matrix computations as in Lemma \ref{vrcvnintpn}, 
it follows that 
\begin{eqnarray*}
\{(-Q \circ \mathcal{S})_{BB}\}^{-1} (b,b') = {\bf{E}}_b \left[ \sum_{i=0}^{\infty} 
\frac{1_{\{S_i < R_A\}} 1_{\{Y_i = b'\}} H_{S_i}}{-q_{b'b'}} \right]. 
\end{eqnarray*}

\noindent
where $R_A$ is the first time (greater than equal to $S_1$) that the Markov chain 
$\{X_t\}_{t \geq 0}$ hits $A$ and $H_{S_i}$ is the ``sign-process" evaluated at the 
$i^{th}$ jump time $S_i$, for $i \geq 0$. Hence, 
\begin{eqnarray*}
{\Sigma}_{BA}^\mathcal{S} {{\Sigma}_{AA}^\mathcal{S}}^{-1} (b,a)
&=& \sum_{b' \in B} {\bf{E}}_b \left[ \sum_{i=0}^{\infty} 1_{\{S_i < R_A\}} 
1_{\{Y_i = b'\}} H_{S_i} \right] \left( \frac{q_{b'a}}{-q_{b'b'}} \right) 
\mathcal{S} (b',a)\\
&=& \sum_{i=0}^{\infty} \sum_{b' \in B} {\bf{E}}_b \left[ 1_{\{S_i < R_A\}} 
1_{\{Y_i = b'\}} 1_{\{Y_{i+1} = a\}} H_{S_{i+1}} \right] 
\end{eqnarray*}

\noindent
The previous equality follows by conditioning, the Markov property and 
$P_{b'} \{ Y_1 = a \} = \frac{q_{b'a}}{-q_{b'b'}} = p(b',a)$. 

\smallskip

\noindent
Hence, 
\begin{eqnarray*}
{\Sigma}_{BA}^\mathcal{S} {{\Sigma}_{AA}^\mathcal{S}}^{-1} (b,a) 
&=& {\bf{E}}_b \left[ \sum_{i=0}^{\infty} 1_{\{S_i < R_A, Y_i \in B, Y_{i+1} = a\}} 
H_{S_{i+1}} \right]\\
&=& {\bf{E}}_b \left[ \sum_{i=0}^{\infty} 1_{\{R_A = S_{i+1}, Y_{i+1} = a\}} 
H_{S_{i+1}} \right] 
\end{eqnarray*}

\noindent
It follows that, 
\begin{equation} \label{coeffidnty}
{\Sigma}_{BA}^\mathcal{S} {{\Sigma}_{AA}^\mathcal{S}}^{-1} (b,a) =  \{(-Q \circ 
\mathcal{S})_{BB}\}^{-1} (-Q \circ \mathcal{S})_{BA} = {\bf{E}}_b \left[ 1_{\{X_{R_A} = 
a\}} H_{R_A} \right]. 
\end{equation}

\noindent
Hence, $\forall b \in B$, 
\begin{equation} \label{dknprssgsn}
{\bf{E}}[Z_b^\mathcal{S} \mid Z_a^\mathcal{S}, \; a \in A] = \sum_{a \in A} {\bf{E}}_b 
\left[ 1_{\{X_{R_A} = a\}} H_{R_A} \right] Z_a^\mathcal{S}. 
\end{equation}

\smallskip

\noindent
Also, 
\begin{eqnarray*}
Var[Z_B^\mathcal{S} \mid Z_A^\mathcal{S}] 
&=& {\Sigma}_{BB}^\mathcal{S} - {\Sigma}_{BA}^\mathcal{S} {{\Sigma}_{AA}^\mathcal{S}}^{-1} 
{\Sigma}_{AB}^\mathcal{S}\\
&=& \{(-Q \circ \mathcal{S})_{BB}\}^{-1} 
\end{eqnarray*}

\noindent
Hence, 
\begin{eqnarray*}
Cov(Z_b^\mathcal{S}, Z_{b'}^\mathcal{S} \mid Z_a^\mathcal{S}, \; a \in A) = {\bf{E}}_b 
\left[ \sum_{i=0}^{\infty} \frac{1_{\{S_i < R_A\}} 1_{\{Y_i = b'\}} H_{S_i}}{-q_{b'b'}} 
\right], \; \forall b,b' \in B. 
\end{eqnarray*}

\smallskip

\noindent
It follows that, 
\begin{equation} \label{dknprssvar}
Cov(Z_b^\mathcal{S}, Z_{b'}^\mathcal{S} \mid Z_a^\mathcal{S}, \; a \in A) = {\bf{E}}_b 
\left[ \int_0^{\infty} 1_{\{X_s = b'\}} H_s 1_{\{s < R_A\}} ds \right], \; \forall b,b' \in 
B. 
\end{equation}

\subsection{The infinite case}
We now deal with the case when $\mathcal{X}$ is countably infinite. To prove 
(\ref{dknimphsgn}) in this case, arbitrarily fix a finite subset $\mathcal{X}_f \subset 
\mathcal{X}$. Note that the variance-covariance matrix for $\{Z_x^\mathcal{S}\}_{x \in 
\mathcal{X}_f}$ is given by 
$$
\Sigma_f := (-(Q \circ \mathcal{S})_{\mathcal{X}_f \mathcal{X}_f} + (Q \circ 
\mathcal{S})_{\mathcal{X}_f \mathcal{X}_f^c} \{(Q \circ \mathcal{S})_{\mathcal{X}_f^c 
\mathcal{X}_f^c}\}^{-1} (Q \circ \mathcal{S})_{\mathcal{X}_f^c \mathcal{X}_f})^{-1}. 
$$

\noindent
Hence for $x,y \in \mathcal{X}_f$, it follows by a similiar calculation leading to 
(\ref{gaussinteg}) that for arbitrary $d_u \geq 0, \; u \in \mathcal{X}_f$, 
$$
{\bf{E}} \left[ Z_x^\mathcal{S} Z_y^\mathcal{S} e^{-\sum_{u \in \mathcal{X}_f} d_u 
\frac{(Z_u^\mathcal{S})^2}{2}} \right] = {\bf{E}} \left[ e^{-\sum_{u \in \mathcal{X}_f} d_u 
\frac{(Z_u^\mathcal{S})^2}{2}} \right] (D_f + \Sigma_f^{-1})^{-1} (x,y). 
$$

\noindent
Here $D_f$ is a diagonal matrix of dimension $|\mathcal{X}_f|$ with diagonal entries $d_u, 
\; u \in \mathcal{X}_f$. We now prove a claim which will help us prove that $(D_f + 
\Sigma_f^{-1})^{-1} (x,y)$ is indeed $(D - Q \circ \mathcal{S})^{-1} (x,y)$, where $D$ is a 
diagonal matrix of dimension $|\mathcal{X}|$ with diagonal entry $d_u$ if $u \in 
\mathcal{X}_f$ and $0$ otherwise. 
\newtheorem{claim}{Claim}[section]
\begin{claim} \label{infstidnty}
Let $A \subset \mathcal{X}$ be finite and 
$$
(Q \circ \mathcal{S})^A := (Q \circ \mathcal{S})_{AA} - (Q \circ \mathcal{S})_{AB} \{(Q 
\circ \mathcal{S})_{BB}\}^{-1} (Q \circ \mathcal{S})_{BA}. 
$$

\noindent
If $Q$ satisfies (\ref{transient1}) or (\ref{transient2}), 
\begin{eqnarray*}
&(a)& \; (Q \circ \mathcal{S})^A (a,a') = \cases{q_{aa}(1 - {\bf{E}}_a [1_{\{X_{R_A} = a\}} 
H_{R_A}]) & \text{if $a = a'$} \cr -q_{aa} {\bf{E}}_a [1_{\{X_{R_A} = a'\}} H_{R_A}] & 
\text{if $a \neq a'$}}\\
&\mbox{and,}& \\
&(b)& \; \{(Q \circ \mathcal{S})^A\}^{-1} (a,a') = (Q \circ \mathcal{S})^{-1} (a,a') \; 
\forall a,a' \in A. 
\end{eqnarray*}
\end{claim}

\noindent
{\it Proof} Throughout the proof, the absolute convergence for various infinite sums will 
be taken care of by the fact that $\sum_{n=0}^\infty P^n < \infty$ (because $Q$ satisifes 
(\ref{transient1}) or (\ref{transient2})). Note that the prediction formulas derived for 
the finite case in Section $4.1$ go through for the infinite case as well, if $A$ is 
finite. Recall that $R_A$ is the first time the process $\{X_t\}_{t \geq 0}$ hits $A$ 
after the initial state. Hence from (\ref{coeffidnty}) and strong Markov property, 
\begin{eqnarray*}
-(Q \circ \mathcal{S})^A (a,a') 
&=& -q_{aa'} - \sum_{b \in B} q_{ab} \mathcal{S} (a,b) {\bf{E}}_b [1_{\{X_{R_A} = a'\}} 
H_{R_A}]\\
&=& -q_{aa'} + q_{aa} {\bf{E}}_a [1_{\{X_{R_A} = a'\}} 1_{\{R_A > S_1\}} H_{R_A}]\\
&=& \cases{-q_{aa}(1 - {\bf{E}}_a [1_{\{X_{R_A} = a\}} H_{R_A}]) & \text{if $a = a'$} \cr 
q_{aa} {\bf{E}}_a [1_{\{X_{R_A} = a'\}} H_{R_A}] & \text{if $a \neq a'$}} 
\end{eqnarray*}

\noindent
This completes the proof of $(a)$. This also gives $(Q \circ \mathcal{S})^A = Q_{AA}^{diag} 
(I - P^A)$ where $P^A (a,a') := {\bf{E}}_a [1_{\{X_{R_A} = a'\}} H_{R_A}], \; \forall a,a' 
\in A$, and $Q^{diag}$ is the diagonal matrix with the same diagonal entries as $Q$. Let 
$R_A^n$ denote the time of $n^{th}$ return to $A$. It follows that, 
\begin{eqnarray*}
\{(Q \circ \mathcal{S})^A\}^{-1} (a,a') 
&=& \frac{1}{q_{aa}} \sum_{n=0}^\infty (P^A)^n (a,a')\\
&=& \frac{1}{q_{aa}} \sum_{n=0}^\infty \sum_{{a_0, a_1, ... , a_n \in A} \atop {a_0 = a, 
a_n = a'}} \prod_{i=0}^{n-1} {\bf{E}}_{a_i} [1_{\{X_{R_A} = a_{i+1}\}} H_{R_A}]\\
&=& \frac{1}{q_{aa}} \sum_{n=0}^\infty {\bf{E}}_a [1_{\{X_{R_A^n} = a'\}} H_{R_A^n}] 
\end{eqnarray*}

\noindent
The previous equality follows by the definition of $\{H_t\}_{t \geq 0}$ and repeated 
application of the strong Markov property. Observing that $X_{S_i} = a'$ only if $S_i = 
R_A^n$ for some $n \geq 1$, we get that, 
\begin{eqnarray*}
\{(Q \circ \mathcal{S})^A\}^{-1} (a,a') 
&=& \frac{1}{q_{aa}} \sum_{i=0}^\infty {\bf{E}}_a [1_{\{X_{S_i} = a'\}} H_{S_i}]\\
&=& -{\bf{E}}_a \left[ \int_0^\infty 1_{\{X_s = a'\}} H_s ds \right]\\
&=& (Q \circ \mathcal{S})^{-1} (a,a') 
\end{eqnarray*}

\noindent
The previous equality follows from Lemma \ref{vrcvnintpn}. The proof of $(b)$ is now 
complete. \hspace{20mm} $\blacksquare$ 

\medskip

\indent
Note that $-(D - Q)$ is a generator matrix that satisfies (\ref{transient1}) or 
(\ref{transient2}) (because $Q$ satisfies one of these conditions). Applying Claim 
\ref{infstidnty} for $-(D - Q)$ with $A = \mathcal{X}_f$, we get that, 
$$
(D_f + \Sigma_f^{-1})^{-1} (x,y) = ((D - Q) \circ \mathcal{S})^{-1} (x,y). 
$$

\noindent
By imitating the proof of Lemma \ref{lpctmidnty} we get 
$$
((D - Q) \circ \mathcal{S})^{-1} (x,y) = \int e^{-\sum_{u \in \mathcal{X}_f} d_u 
l_\infty^u} H_\infty d \mu_{xy}. 
$$

\noindent
Combining everything, 
$$
{\bf{E}} \left[ Z_x^\mathcal{S} Z_y^\mathcal{S} e^{-\sum_{u \in \mathcal{X}_f} d_u 
\frac{(Z_u^\mathcal{S})^2}{2}} \right] = \int {\bf{E}} \left[ e^{-\sum_{u \in 
\mathcal{X}_f} d_u \left( \frac{(Z_u^\mathcal{S})^2}{2} + l_\infty^u \right)} \right] d 
\mu_{xy}. 
$$

\noindent
Note that the set of functions $F_{\underline{d}} (\underline{w}) = e^{-\sum_{u \in 
\mathcal{X}} d_u w_u}$, where $\underline{d} = \{d_u\}_{u \in \mathcal{X}}$ is arbitrary 
with $d_u > 0$ for finitely many $u \in \mathcal{X}$ and $d_u = 0$ otherwise, generate the 
Borel $\sigma$-algebra in $\mathbb{R}^{|\mathcal{X}|}$ and they form a closed class under 
multiplication. Also, the set of functions $F$ for which 
$$
{\bf{E}} \left[ Z_x^\mathcal{S} Z_y^\mathcal{S} F \left( \frac{(Z_u^\mathcal{S})^2}{2}, \; 
u \in \mathcal{X} \right) \right] = \int {\bf{E}} \left[ F \left( 
\frac{(Z_u^\mathcal{S})^2}{2} + l_\infty^u, \; u \in \mathcal{X} \right) H_\infty \right] d 
\mu_{xy}, 
$$

\noindent
is a linear space closed under bounded convergence and under monotone convergence. Hence 
for each bounded Borel measurable function $F: \mathbb{R}^{|\mathcal{X}|} \rightarrow 
\mathbb{R}$, 
$$
{\bf{E}} \left[ Z_x^\mathcal{S} Z_y^\mathcal{S} F \left( \frac{(Z_u^\mathcal{S})^2}{2}, \; 
u \in \mathcal{X} \right) \right] = \int {\bf{E}} \left[ F \left( 
\frac{(Z_u^\mathcal{S})^2}{2} + l_\infty^u, \; u \in \mathcal{X} \right) H_\infty \right] d 
\mu_{xy}. 
$$

\subsection{Conditional Independence Property}
There is another interesting property of Dynkin's isomorphism which is preserved after 
introducing a ``sign" matrix $\mathcal{S}$. Let $Q$ be the generator of a 
continuous time Markov process $\{X_t\}_{t \geq 0}$, with a countable state space 
$\mathcal{X}$. Assume $Q^{-1}$ exists and $Q$ is symmetric. Let $\{Z_x^\mathcal{S}\}_{x \in 
\mathcal{X}}$ be a zero mean Gaussian field with variance-covariance matrix 
$\Sigma^\mathcal{S} = (-Q \circ \mathcal{S})^{-1}$. 
\begin{lem}
Let $A,B,C$ be disjoint subsets of the state space $\mathcal{X}$, such that to go from 
any state in $A$ to any state in $C$, the Markov process $\{X_t\}_{t \geq 0}$ has to 
pass through $B$. Then conditioned on $Z_B^\mathcal{S} := \{Z_b^\mathcal{S}\}_{b \in B}$, 
the Gaussian random vectors $Z_A^\mathcal{S} := \{Z_a^\mathcal{S}\}_{a \in A}$ and 
$Z_C^\mathcal{S} := \{Z_c^\mathcal{S}\}_{c \in C}$ are independent. 
\end{lem}

\noindent
{\it Proof} Fix $a \in A$ and $c \in C$ arbitrarily. By (\ref{dknprssvar}), 
\begin{eqnarray*}
Cov(Z_a^\mathcal{S}, Z_c^\mathcal{S} \mid Z_b^\mathcal{S}, \; b \in B) = {\bf{E}}_a \left[ 
\int_0^{\infty} 1_{\{X_s = c\}} 1_{\{s < R_B\}} H_s ds \right]. 
\end{eqnarray*}

\noindent
Since the Markov process $\{X_t\}_{t \geq 0}$ has to pass through the set $B$ to go 
from the state $a$ to the state $c$, 
$$
1_{\{X_s = c\}} 1_{\{s < R_B\}} = 0 \mbox{ under } {\bf{P}}_a. 
$$

\noindent
Hence, 
\begin{eqnarray*}
Cov(Z_a^\mathcal{S}, Z_c^\mathcal{S} \mid Z_b^\mathcal{S}, \; b \in B) = 0. 
\end{eqnarray*}

\noindent
Since $a \in A$ and $c \in C$ were arbitrarily fixed, it follows that $Z_A^\mathcal{S}$ and 
$Z_C^\mathcal{S}$ are uncorrelated given $Z_B^\mathcal{S}$. Two random vectors having a 
joint Gaussian distribution are independent iff they are uncorrelated. Hence, 
$Z_A^\mathcal{S}$ and $Z_C^\mathcal{S}$ are independent given $Z_B^\mathcal{S}$. 
\hspace{14mm} $\blacksquare$ 

\smallskip

\noindent
Bolthausen \cite{bolthausen} uses this property in his analysis of the 
Gaussian free field. 

\subsection{An Example: Ornstein-Uhlenbeck Process on $\mathbb{N}$}

\noindent
Consider the Ornstein-Uhlenbeck process $\{Z_i\}_{i \in \mathbb{N}}$ defined by 
$$
Z_1 = {\varepsilon}_1, \; \; Z_i = aZ_{i-1} + {\varepsilon}_i \; \forall i \geq 
2, \mbox{ where } \{{\varepsilon}_i\}_{i \geq 1} \mbox{ are i.i.d. } N(0,1). 
$$

\noindent
Let $\Sigma$ denote the variance-covariance matrix of $\{Z_i\}_{i \in \mathbb{N}}$. 
Then, 
\begin{equation} \label{dfnvacvmtx}
\Sigma (k,l) = \cases{a^{l-k} \sum_{i=1}^k a^{2(k-i)} & \text{if $k \leq l$}, \cr 
\Sigma (l,k) & \text{if $l < k$}.} 
\end{equation}

\noindent
After some manipulations, we can establish that $\Sigma = -Q^{-1}$ where 
\begin{equation} \label{dfngntrmtx}
Q(k,l) = \cases{-(1 + a^2) & \text{if $k = l$}, \cr 
a & \text{if $k = l \pm 1$}, \cr 
0 & \text{otherwise}.} 
\end{equation}

\noindent
If $a > 0$, then $Q$ is the generator of a birth and death process. Hence, the 
Ornstein-Uhlenbeck process is connected to the birth and death process with 
generator $Q$ by Dynkin's isomorphism. 

\noindent
Suppose we introduce a ``sign" matrix $\mathcal{S}$ (as described earlier in 
this section) and for $a > 0$, look at a Gaussian field 
$\{Z_i'\}_{i \in \mathbb{N}}$ defined by 
$$
Z_1' = {\varepsilon}_1, \; \; Z_i' = \mathcal{S} (i-1,i) aZ_{i-1}' + 
{\varepsilon}_i, \; \forall i \geq 2, \mbox{ where } 
\{{\varepsilon}_i\}_{i \geq 1} \mbox{ are i.i.d. } N(0,1). 
$$

\noindent
It follows after some manipulations that the variance-covariance matrix of 
$\{Z_i'\}_{i \in \mathbb{N}}$ is given by $(-Q \circ \mathcal{S})^{-1}$, 
where $Q$ is as specified in (\ref{dfngntrmtx}). 

\noindent
If $\mathcal{S} (k,l) = -1 \; \forall k \neq l$, then 
$$
Z_1' = {\varepsilon}_1, \; \; Z_i' = -a Z_{i-1}' + {\varepsilon}_i, \; 
\forall i \geq 2. 
$$

\noindent
Hence, if $a > 0$, then the Ornstein-Uhlenbeck process on $\mathbb{N}$ with 
parameter $-a$ is associated to a birth and death process with generator 
$Q$ in (\ref{dfngntrmtx}), by Dynkin's isomorphism with ``sign" matrix 
$\mathcal{S}$ such that $\mathcal{S} (k,l) = -1 \; \forall k \neq l$. 

\section{An Algorithm for Computing the Prediction Coefficients}
We present an algorithm for computing the prediction coefficients ${\bf{E}}_b 
\left[ 1_{\{X_{R_A} = a\}} H_{R_A} \right]$ in (\ref{dknprssgsn}) for calculating 
${\bf{E}}[Z_b^\mathcal{S}\mid Z_a^\mathcal{S}, \; a \in A]$. This algorithm can be 
described in two ways: 
\begin{itemize}
\item Graph theoretic description. 

\noindent
Construct a graph on the vertex set $\mathcal{X}$ by putting an edge of weight 
$Q(x,y) \mathcal{S} (x,y)$ between vertices $x$ and $y, \; \forall x \neq y \in 
\mathcal{X}$. If any of these weights are $0$, that by default means no edge is 
put between the corresponding vertices. Put a loop of weight $-Q(x,x)$ at each 
vertex $x \in \mathcal{X}$. Since we want to predict $Z_b^\mathcal{S}$ given 
$\{Z_a^\mathcal{S}\}_{a \in A}$ we now proceed to remove all vertices not in $\{b\} \cup A$ 
from this graph in a sequential fashion. Choose any vertex, say $z$ not in 
$\{b\} \cup A$. If we remove $z$, i.e. if we behave as if $z$ does not exist in 
the graph, this leads to forming a new edge between every pair $x$ and $y$ such 
that $x$ and $z$, as well as $y$ and $z$ share an edge. The weight of this 
new edge is the product of the weights of these two edges divided by the 
weight of the loop at $z$. If there is already an edge between $x$ and $y$, add 
the weight of this new edge to the existing one and combine them into one edge. 
Perform this procedure with all $x$ and $y$ sharing an edge with $z$ (including 
the case $x = y$). So we get a new graph with vertex set $\mathcal{X} \setminus 
\{z\}$ and edge set as described above. Note that $-Q \circ \mathcal{S}$ is a 
diagonally dominant matrix with positive diagonal entries and hence for the 
old graph, the weight of the loop at any vertex dominates the sum of the 
absolute weight values of the the edges emanating from that vertex. As we will 
see later, the new grpah has the same property. We continue choosing vertices 
and removing them by using the above procedure until we are left with the 
vertex set $A \cup \{b\}$. The coefficient of $Z_a^\mathcal{S}$ in 
${\bf{E}}[Z_b^\mathcal{S} \mid Z_a^\mathcal{S}, \; a \in A]$ is precisely the weight of the 
edge joining $a$ and $b$ divided by the weight of the loop at $b$, for every $a$ in A. 
\item Analytic description. 

\noindent
We can describe the above algorithm analytically as follows: 
\begin{enumerate}
\item Start with $V = \mathcal{X}, \; M = -Q \circ \mathcal{S}$. 
\item Choose $z \in V \setminus \{A \cup \{b\}\}$. 
\item $M(x,y) = M(x,y) - \frac{M(x,z)M(y,z)}{M(z,z)} \; \forall x,y \in V$. 
\item Remove the $z^{th}$ row and the $z^{th}$ column of $M$. 
\item $V \rightarrow V \setminus \{z\}$. If $V \neq A \cup \{b\}$ goto 
step $2$, otherwise stop. The coefficient of $Z_a^\mathcal{S}$ in ${\bf{E}}[Z_b^\mathcal{S} 
\mid Z_a^\mathcal{S}, \; a \in A]$ is $- \frac{M(a,b)}{M(b,b)}$ for every $a \in A$. 
\end{enumerate}
\end{itemize}

\noindent
The above description tells us that our algorithm is essentially the sequential 
process of evaluating the Schur complement $(-Q \circ \mathcal{S})_{VV} - 
(-Q \circ \mathcal{S})_{VV^c} \{(-Q \circ \mathcal{S})_{V^c V^c} \}^{-1} 
(-Q \circ \mathcal{S})_{V^c V}$ (and ending at $V = A \cup \{b\}$) by reducing 
rows and columns. Since the Schur complement of a diagonally dominant matrix is 
also diagonally dominant, the matrix $M$ is a diagonally dominant matrix at 
every step of the algorithm. The proof of this algorithm can be obtained 
immediately by observing two facts. Firstly, $\Sigma^\mathcal{S} = (-Q \circ 
\mathcal{S})^{-1}$ implies that 
$$
(\Sigma_{VV}^\mathcal{S})^{-1} = (-Q \circ \mathcal{S})_{VV} - 
(-Q \circ \mathcal{S})_{VV^c} \{(-Q \circ \mathcal{S})_{V^c V^c} \}^{-1} 
(-Q \circ \mathcal{S})_{V^c V} \mbox{ for every } V \subseteq \mathcal{X}. 
$$

\noindent
Hence, when we stop the algorithm, the matrix $M$ is the same as 
$(\Sigma_{VV}^\mathcal{S})^{-1}$ with $V = A \cup \{b\}$. Secondly, if $\underline{Y} \sim 
MVN_n (\underline{0}, \Gamma)$, then 
\begin{equation} \label{invfrmlgsn}
{\bf{E}}[Y_i \mid Y_j, \; j \neq i] = \sum_{j \neq i} \frac{- {\Gamma}^{-1} 
(i,j)}{{\Gamma}^{-1} (i,i)} Y_j. 
\end{equation}

\noindent
Since $\{Z_v^\mathcal{S}\}_{v \in V}$ is $MVN_{|V|} (\underline{0}, 
{\Sigma}_{VV}^\mathcal{S})$ and 
$M = ({\Sigma}_{VV}^\mathcal{S})^{-1}$, it follows by (\ref{invfrmlgsn}) that 
\begin{eqnarray*}
{\bf{E}}[Z_b^\mathcal{S} \mid Z_a^\mathcal{S}, \; a \in V \setminus \{ b\}] = \sum_{a \in V 
\setminus \{ b\}} \frac{-M(a,b)}{M(b,b)} Z_a^\mathcal{S}. 
\end{eqnarray*}

\noindent
Hence this algorithm is not all that mysterious. If $|\mathcal{X}| = n$, the 
worst case running time of this algorithm is $O(n^3)$. One nice property of 
this algorithm is that at any step of the algorithm with vertex set $V$ and 
corresponding matrix $M$, 
\begin{eqnarray*}
{\bf{E}}[Z_v^\mathcal{S} \mid Z_w^\mathcal{S}, \; w \in V \setminus \{ v\} ] = \sum_{w \in 
V \setminus \{ v\}} \frac{-M(v,w)}{M(v,v)} Z_w^\mathcal{S}. 
\end{eqnarray*}

\noindent
Hence, if we want we can obtain the prediction coefficients given 
$\{Z_v^\mathcal{S}\}_{v \in V}$ for every $V \subseteq \mathcal{X}$ that comes 
up in the course of this algorithm. 

\subsection{An Example of Prediction with Independent Errors at the Observed Values}
Consider the Ornstein-Uhlenbeck process on $\mathbb{N}$ with $a = 1$, i.e. 
$$
Z_1 = {\varepsilon}_1, \; \; Z_i = Z_{i-1} + {\varepsilon}_i \; \forall i \geq 2, 
\mbox{ where } \{{\varepsilon}_i\}_{i \geq 1} \mbox{ are i.i.d. } N(0,1). 
$$

\noindent
This process is same as the Gaussian free field on $\mathbb{N}$. It follows that 
the variance-covariance matrix $\Sigma$ of $\{Z_i\}_{i \in \mathbb{N}}$ is given by 
$$
\Sigma (k,l) = k \wedge l \; \forall k,l \in \mathbb{N}. 
$$

\noindent
Suppose we observe the values of the process in the set $V = \{n_1, n_2, ... , 
n_k\}$ (where $n_i < n_j$ if $i < j$), but with an independent additive error  
${\tilde{\varepsilon}}_i$ at the point $n_i \; i = 1,2, ... ,k$, where  
$\{{\tilde{\varepsilon}}_i\}_{1 \leq i \leq k}$ are i.i.d. $N(0, {\sigma}^2)$. 
With these observations, we want to predict the process 
$\{Z_i\}_{i \in \mathbb{N}}$ i.e. we want to compute the expectation 
$$
{\bf{E}}[Z_n \mid Z_{n_1} + {\tilde{\varepsilon}}_1, Z_{n_2} + 
{\tilde{\varepsilon}}_2, ... , Z_{n_k} + {\tilde{\varepsilon}}_k], \; \forall n \in 
\mathbb{N}. 
$$

\noindent
It is known that 
$$
{\bf{E}}[Z_n \mid Z_{n_i} + {\tilde{\varepsilon}}_i, \; i = 1,2, ... ,k] = 
{\Sigma}_{nV} ({\Sigma}_{VV} + {\sigma}^2 I_{|V|})^{-1} Z_V. 
$$

\noindent
Since $\Sigma(n_i,n) = n_i \wedge n \; \forall i = 1,2, ... ,k$, we would like to 
compute a simplified expression for 
$({\Sigma}_{VV} + {\sigma}^2 I_{|V|})^{-1}$. We utilize the structure of 
${\Sigma}_{VV}$ for this purpose. 
$$
{\Sigma}_{VV} = UDU^T, 
$$

\noindent
where, $U(i,j) = 1_{\{i \geq j\}}, \; \forall 1 \leq i,j \leq k$ and $D$ is a 
diagonal matrix with $D(i,i) = n_i - n_{i-1} \; \forall 1 \leq i \leq k$ (where 
$n_0 = 0$). It follows that 
$$
{\Sigma}_{VV} + {\sigma}^2 I_{|V|} = U \Lambda U^T, 
$$

\noindent
where $\Lambda$ is the tridiagonal matrix with 
\begin{equation}
\Lambda(i,j) = \cases{n_i - n_{i-1} + i{\sigma}^2 & \text{if $i = j$}, \cr 
-{\sigma}^2 & \text{if $|i-j| = 1$}, \cr 
0 & \text{otherwise}.} 
\end{equation}

\noindent
Note that, 
\begin{equation}
U^{-1} (i,j) = \cases{1 & \text{if $i = j$}, \cr 
-1 & \text{if $i = j + 1$}, \cr 
0 & \text{otherwise}.} 
\end{equation}

\noindent
Let $r_0 = {\sigma}^2, \; r_1 = (n_1 + {\sigma}^2) {\sigma}^2, \; r_i = 
(n_i - n_{i-1} + i{\sigma}^2)r_{i-1} - r_{i-2} \mbox{ for } i = 2,3, ... ,k$. 
By the explicit formula for the inverse of a symmetric tridiagonal matrix in 
\cite{schlegeltm}, 
\begin{equation}
{\Lambda}^{-1} (i,j) = \cases{\frac{r_{i-1} r_{k-j}}{r_k} & \text{if $i \leq j$}, \cr 
{\Lambda}^{-1} (j,i) & \text{if $i > j$}.} 
\end{equation}

\noindent
Hence we obtain 
\begin{eqnarray*}
{\bf{E}}[Z_n \mid Z_{n_i} + {\tilde{\varepsilon}}_i, \; i = 1,2, ... ,k] = 
\sum_{i=1}^k (\sum_{j=1}^k {\gamma}_{ij} (n \wedge n_j)) (Z_{n_i} + 
{\tilde{\varepsilon}}_i), 
\end{eqnarray*}

\noindent
where, 
\begin{equation}
{\gamma}_{ij} = \cases{\frac{(r_i - r_{i-1})(r_{k-j} - r_{k-j-1})}{r_k} & \text{if 
$i < j$}, \cr \frac{r_{i-1} (r_{k-i} - 2r_{k-i-1}) + r_i r_{k-i-1}}{r_k} & \text{if 
$i = j$}, \cr {\gamma}_{ji} & \text{if $i > j$}.} 
\end{equation}

\noindent
As is clear from this example, introducing errors leads to non-trivial changes in 
the prediction coefficients. It is hard to find a general formula which expresses 
these changed coefficients in terms of the associated Markov chain.

\end{document}